\documentclass[letterpaper, 10 pt, conference]{ieeeconf}

\usepackage{amsmath,amssymb,amsfonts,mathrsfs, mathtools, bbm, dsfont, mathrsfs,epstopdf}
\usepackage{graphicx}

\let\labelindent\relax

\IEEEoverridecommandlockouts

\usepackage{subcaption}
\usepackage{ifthen}
\usepackage{multirow, multicol}
\usepackage{float}
\usepackage{subeqnarray, array}
\usepackage{enumitem}
\usepackage{xcolor}
\usepackage[font=small]{caption}

\usepackage[hyphens]{url}
\usepackage[breaklinks, colorlinks, linkcolor=purple, anchorcolor=purple, citecolor=purple, urlcolor=purple]{hyperref}
\usepackage{cite}

\usepackage[ruled]{algorithm2e}

\usepackage{tikz}

\DeclareMathAlphabet{\mathpzc}{OT1}{pzc}{m}{it}

\newcommand{\subparagraph}{}
\usepackage[explicit]{titlesec}
\usepackage{ulem}

\titleformat{\subsubsection}
  {\normalfont\itshape}{\thesubsection}{1em}{$\bullet$ {#1}.}
\titlespacing*{\subsubsection}{0pt}{2ex}{2ex}

\graphicspath{{Figures/}}

\newboolean{showcomments}
\setboolean{showcomments}{true}

\newcommand{\bose}[1]{\ifthenelse{\boolean{showcomments}}
{ \textcolor{red}{(Bose says:  #1)}}{}}

\newcommand{\TD}[1]{\ifthenelse{\boolean{showcomments}}
{ \textcolor{red}{(Thinh says:  #1)}}{}}

\newcommand{\redline}{\raisebox{2pt}{\tikz{\draw[-,red,solid,line width = 1.5pt](0,0) -- (6mm,0);}}}
\newcommand{\blueline}{\raisebox{2pt}{\tikz{\draw[-,blue,solid,line width = 1.5pt](0,0) -- (6mm,0);}}}

\graphicspath{{Figures/}}

\newcommand{\beq}{\begin{equation}}
\newcommand{\eeq}{\end{equation}}
\newcommand{\beqa}{\begin{eqnarray}}
\newcommand{\eeqa}{\end{eqnarray}}
\newcommand{\beqan}{\begin{eqnarray*}}
\newcommand{\eeqan}{\end{eqnarray*}}

\newcommand{\vnorm}[1]{\left\|#1\right\|}

\newcommand\T{{\mathpalette\raiseT\intercal}}
\newcommand\raiseT[2]{\raisebox{0.25ex}{$#1#2$}
}

\newcommand{\Dset}{\mathbb{D}}

\newcommand{\Hset}{\mathbb{H}}

\newcommand{\Rset}{\mathbb{R}}

\newcommand{\Xset}{\mathbb{X}}

\newcommand{\Ocal}{{\cal O}}

\newcommand{\Gfrak}{\mathfrak{G}}

\newcommand{\argmin}{\mathop{\rm argmin}}

\newcommand{\bone}{\mathbf{1}}

\renewcommand{\v}[1]{{\mathbf{#1}}}

\newcommand{\ol}[1]{\ensuremath{\overline{{#1}}}}

\newcommand{\sgn}{\text{sgn }}

\newcounter{l1}
\newcounter{l2}
\newcounter{l3}
\newcommand{\bdotlist}{\begin{list}{$\bullet$}{}}
\newcommand{\bboxlist}{\begin{list}{$\Box$}{}}
\newcommand{\bbboxlist}{\begin{list}{\raisebox{.005in}{{\tiny
$\blacksquare$ \ \ }}}{}}
\newcommand{\bdashlist}{\begin{list}{$-$}{} }
\newcommand{\blist}{\begin{list}{}{} }
\newcommand{\barablist}{\begin{list}{\arabic{l1}}{\usecounter{l1}}}
\newcommand{\balphlist}{\begin{list}{(\alph{l2})}{\usecounter{l2}}}
\newcommand{\bAlphlist}{\begin{list}{\Alph{l2}.}{\usecounter{l2}}}
\newcommand{\bdiamlist}{\begin{list}{$\diamond$}{}}
\newcommand{\bromalist}{\begin{list}{(\roman{l3})}{\usecounter{l3}}}

\newtheorem{theorem}{Theorem}

\newtheorem{remark}{Remark}

\renewcommand{\bone}{\mathds{1}}

\title{\textbf{\LARGE Convergence of the Iterates in Mirror Descent Methods}}

\author{Thinh T. Doan  \qquad Subhonmesh Bose \qquad D. Hoa Nguyen \qquad Carolyn L. Beck\thanks{S.Bose and T. T. Doan are with the Dept. of Electrical and Computer Engg., and C. L. Beck is with the Dept. of Industrial and Systems Engg. at the University of Illinois at Urbana-Champaign, Urbana, IL, USA. D. H. Nguyen is with the International Institute of Carbon-Neutral Energy Research and the Institute of Mathematics for Industry at Kyushu University, Japan. Emails: \{ttdoan2, boses, beck3\}@illinois.edu, hoa.nd@i2cner.kyushu-u.ac.jp.}}

\begin{document}

\maketitle

\begin{abstract}
We consider centralized and distributed mirror descent algorithms over a finite-dimensional Hilbert space, and prove that the problem variables converge to an optimizer of a possibly nonsmooth function when the step sizes are square summable but not summable. Prior literature has focused on the convergence of the {\em function value} to its optimum. However, applications from distributed optimization and learning in games require the convergence of the variables to an optimizer, which is generally not guaranteed without assuming strong convexity of the objective function. We provide numerical simulations comparing entropic mirror descent and standard subgradient methods for the robust regression problem.
\end{abstract}


\section{Introduction}

The method of Mirror Descent ({\sf MD}), originally proposed by Nemirovski and Yudin \cite{Nemirovsky1983}, is a primal-dual method for solving constrained convex optimization problems.  {\sf MD} is fundamentally a subgradient projection ({\sf SGP}) algorithm that allows one to exploit the geometry of an optimization problem 
through an appropriate choice of a strongly convex function \cite{BeckT2003}.  This method not only generalizes the standard gradient descent ({\sf GD}) method, but also achieves a better convergence rate.  In addition, 
{\sf MD} is applicable to optimization problems in Banach spaces where {\sf GD} is not \cite{Bubeck2015}.  

Of more recent interest, {\sf MD} has been shown to be useful for efficiently solving large-scale optimization problems.   In general, SGP algorithms are simple to implement, however they are typically slow to converge due to the fact that they are based in Euclidean spaces, and through a projection operator are inevitably tied to the geometry of these spaces.  As a result, their convergence rate may be directly tied to the dimension $d$ of the underlying Euclidean space in which the problem variables reside.  Alternatively, {\sf MD} can be adapted, or more specifically tailored to the geometry
of the underlying problem space, potentially allowing for an improved rate of convergence; see \cite{BenTalMN2001} for an early example.  Because of these notable potential benefits, {\sf MD} has experienced significant recent attention for applications to large-scale optimization and machine learning problems in both the continuous  and discrete time settings \cite{Raginsky2012, Krichene2015}, both the deterministic and stochastic scenarios \cite{Duchi2011,Nedic2014,Rabbat2015,Nokleby2017}, and both the centralized and distributed contexts \cite{Rabbat2015,JLi2016}. 
{\sf MD} has also been applied to a variety of practical problems, e.g., game-theoretic applications \cite{Zhou2017a},  and multi-agent distributed learning problems \cite{Hall2015,Ledva2015,Nedic2016,Shahrampour2017,Zhou2017b}.

Hitherto, most of these prior studies have focused on studying the convergence rate of {\sf MD}. In particular, if the step sizes are properly selected then {\sf MD} can achieve a convergence rate of $\mathcal{O}(1/k)$ or $\mathcal{O}(1/\sqrt{k})$ for strongly convex or convex objective functions, respectively, 
\cite{Nemirovsky1983,Nedic2014}. However, the convergence of the objective function value does not, in general, imply the convergence of the sequence of variables to an optimizer\footnote{One can only show the convergence of the sequence of these variables to the optimal set when this set is bounded.}.  To the best of the authors' knowledge, there has not been any prior work establishing the convergence of these variables to an optimizer.  Our motivation for pursuing a study of the convergence to an optimizer arises from potential applications in Distributed Lagrangian ({\sf DL}) methods and Game Theory.  Specifically, in the context of {\sf DL} methods, we can apply distributed subgradient methods, or preferably distributed {\sf MD} methods, to find the solution to the dual problem.  In this setting, convergence to the dual optimizer is needed to complete the convergence analysis of {\sf DL} methods \cite{DoanB2017, DoanBB2017}. To motivate our study from a game theoretic viewpoint, note that the dynamics of certain natural learning strategies in routing games have been identified as the dynamics of centralized mirror descent in the strategy space of the players; see \cite{KricheneKB2015} for an example. In that context, convergence of the learning dynamics to the Nash equilibria (the minimizers of a convex potential function of the routing game) is critical; convergence to the optimal function value is not enough.
 

In this paper, our main contribution is thus a proof of convergence to an optimizer in the {\sf MD} method, where the objective function is convex and not necessarily differentiable; we consider both the centralized and distributed settings.



\section{Centralized Mirror Descent}
\label{sec:centralized}

Let $\left(\Hset, \langle \cdot, \cdot \rangle \right)$ describe a finite-dimensional Hilbert space over reals, and $\Xset$ be a closed convex subset of $\Hset$. Consider a possibly nonsmooth convex and continuous function $f:\Xset \to \Rset$ that we seek to minimize via {\sf MD}, starting from $\v{x}_0 \in \Xset$. We assume throughout that $f$ is finite-valued, and its effective domain contains $\Xset$. Further, we assume throughout that $f$ has at least one finite optimizer over $\Xset$. 


To precisely define {\sf MD}, consider a continuously differentiable $\mu$-strongly convex function on an open convex set $\Dset$ whose the closure contains $\Xset$. By that, we mean $\psi$ satisfies
$$  \psi(\v{y}) \geq \psi(\v{x}) + \left\langle \nabla\psi(\v{x}), \v{y} - \v{x} \right\rangle + \frac{\mu}{2} \vnorm{\v{y}-\v{x}}^2,\;\; \forall\; \v{x},\v{y}\in\Dset. 
$$
Here, $\nabla\psi$ denotes the gradient of $\psi$ which is assumed to be diverged on the boundary of $\Dset$. In addition, $\vnorm{\cdot}$ is the norm induced by the inner product. Define the \emph{Bregman divergence} associated with $\psi$ for all $\v{x}$ and $\v{y}$ in $\Xset$ as
\begin{align}
D_\psi(\v{y},\v{x}) = \psi(\v{y}) - \psi(\v{x}) - \langle \nabla\psi(\v{x}), \v{y} - \v{x} \rangle. \label{prob:bregman}
\end{align}

Equipped with this notation, {\sf MD} prescribes the following iterative dynamics, starting from some $\v{x}_0\in\Xset$.
\begin{align}
\v{x}_{k+1} = \argmin_{\v{z}\in\Xset}\left\{ \langle \nabla f(\v{x}_k), \v{z} - \v{x}_k \rangle+ \frac{1}{\alpha_k}D_\psi(\v{z},\v{x}_k) \right\}.\label{alg:MirrorDescent} 
\end{align}
Here, $\nabla f(\v{x}_k) \in \partial f(\v{x}_k)$ is an arbitrary subgradient of $f$ at $\v{x}_k$, the collection of which comprise the subdifferential set $\partial f(\v{x}_k)$, defined as 
$$ \partial f(\v{x}_k) := \{ \v{g} \in \Hset \ \vert \ f(\v{y}) \geq f(\v{x}_k) + \langle \v{g}, \v{y} - \v{x}_k \rangle 
\ \text{for all } \v{y} \in \Xset \}.$$ 


We note that {\sf MD} enjoys an optimal $\Ocal\left(\frac{1}{\sqrt{k}}\right)$ convergence rate for nonsmooth functions \cite{Nemirovsky1983,Bubeck2015}, i.e.,
$$ f\left(\frac{1}{k}\sum_{t=1}^{k}\v{x}_t\right) - f^* \leq \frac{C}{\sqrt{k}}, $$
where $f^*$ is the optimal value of $f$ over $\Xset$. Of interest to us in this work is the possible convergence of the problem variables themselves, i.e., whether $\v{x}_k$ converges to an optimizer for a suitable choice of step sizes $\alpha_k$. In the remainder of this section, we prove such a convergence result for centralized ${\sf MD}$, and extend this to a distributed setting in the next section.

\begin{theorem}
\label{thm:CMD}
Suppose 
\begin{itemize}
\item $f$ is $L$-Lipschitz continuous over $\Xset$, and 
\item $\{\alpha_k\}_{k=0}^\infty$ defines a nonincreasing sequence of positive step sizes that is square-summable, but not summable, i.e.,
$\sum_{k=0}^{\infty} \alpha_k = \infty,$ $\sum_{k=0}^{\infty} \alpha_k^2 < \infty$.
\end{itemize}
Then, $\lim_{k \to \infty} \v{x}_k$ optimizes $f$ over $\Xset$ for $\v{x}_k$'s generated by {\sf MD} in \eqref{alg:MirrorDescent}.
\end{theorem}

In proving the result, the following two properties of Bregman divergence will be useful. Their proofs are straightforward from its definition in \eqref{prob:bregman}.
\begin{align}
\begin{aligned}
& D_\psi(\v{y}, \v{x}) - D_\psi(\v{y}, \v{z}) - D_\psi(\v{z}, \v{x}) \\
& \qquad =  \left\langle \nabla\psi(\v{z})-\nabla\psi(\v{x}), \v{y}- \v{z} \right\rangle, \\
& D_\psi(\v{z}, \v{x}) \geq \frac{\mu}{2}\vnorm{\v{z}-\v{x}}^2
\end{aligned}
\label{eq:prop.bregman}
\end{align}
for arbitrary $\v{x}$, $\v{y}$, $\v{z}$ in $\Xset$.

%
%
%



\subsection{Proof of Theorem \ref{thm:CMD}}
Our proof proceeds in two steps. We first show that consecutive iterates satisfy
\begin{align}
& D_\psi( \v{z},\v{x}_{k+1} ) - D_\psi( \v{z},\v{x}_k) \notag \\
& \qquad \qquad \leq \alpha_k \langle \nabla f(\v{x}_k), \v{z}-\v{x}_k\rangle + \frac{\alpha_k^2 L^2}{2\mu}\label{lemma1:MainIneq}
\end{align}
for each $\v{z}\in\Xset$. We then deduce the result from \eqref{lemma1:MainIneq}.

\subsubsection*{Proof of \eqref{lemma1:MainIneq}} 
The optimality of $\v{x}_{k+1}$ in \eqref{alg:MirrorDescent} implies
\begin{align}
\langle \alpha_k \nabla f(\v{x}_k) + \nabla_1 D_\psi(\v{x}_{k+1},\v{x}_k), \v{z}-\v{x}_{k+1} \rangle \geq 0.
\label{lemma1:Eq1}
\end{align}
Here, $\nabla_1 D_\psi$ stands for the derivative of the Bregman divergence with respect to the first coordinate. 
The properties of the divergence in \eqref{eq:prop.bregman} yield
\begin{align*}
& \langle \nabla_1 D_\psi( \v{x}_{k+1},\v{x}_k), \v{z} - \v{x}_{k+1} \rangle \\ 
& \qquad = \langle \nabla \psi(\v{x}_{k+1}) - \nabla \psi(\v{x}_k), \v{z} - \v{x}_{k+1} \rangle \nonumber\\
& \qquad = D_\psi(\v{z},\v{x}_{k}) - D_\psi( \v{z},\v{x}_{k+1}) - D_\psi(\v{x}_{k}, \v{x}_{k+1})\nonumber\\
& \qquad \leq D_\psi(\v{z},\v{x}_{k} ) - D_\psi(\v{z},\v{x}_{k+1}) - \frac{\mu}{2}\vnorm{\v{x}_{k+1}-\v{x}_{k}}^2.
\end{align*}
Substituting the above relation in \eqref{lemma1:Eq1}, we get
\begin{align}
& D_\psi(\v{z},\v{x}_{k+1} ) - D_\psi(\v{z},\v{x}_{k} ) \notag \\
& \qquad \leq \alpha_k \langle \nabla f(\v{x}_k), \v{z}-\v{x}_{k+1} \rangle  - \frac{\mu}{2}\vnorm{\v{x}_{k+1}-\v{x}_{k}}^2 \notag \\
& \qquad = \alpha_k \langle \nabla f(\v{x}_k), \v{z}-\v{x}_{k} \rangle - \frac{\mu}{2}\vnorm{\v{x}_{k+1}-\v{x}_{k}}^2 \notag \\ 
& \qquad \qquad + \alpha_k \langle \nabla f(\v{x}_k), \v{x}_k-\v{x}_{k+1} \rangle. \label{eq:ineq.2}
\end{align}
An appeal to Cauchy-Schwartz and arithmetic-geometric mean inequalities allows us to bound the last term on the right hand side of the above inequality, as follows,
\begin{align*}
\langle \alpha_k \nabla f(\v{x}_k), \v{x}_k-\v{x}_{k+1} \rangle \leq \frac{\alpha_k^2}{2\mu} \vnorm{\nabla f(\v{x}_k)}^2 + \frac{\mu}{2} \vnorm{\v{x}_{k+1}-\v{x}_k}^2.
\end{align*}
The above inequality and \eqref{eq:ineq.2} together imply
\begin{align*}
& D_\psi(\v{z},\v{x}_{k+1} ) - D_\psi(\v{z},\v{x}_{k} ) \notag \\
& \qquad \leq \alpha_k \langle \nabla f(\v{x}_k), \v{z}-\v{x}_{k} \rangle +  \frac{\alpha_k^2}{2\mu} \vnorm{\nabla f(\v{x}_k)}^2.
\end{align*}
Lipschitz continuity of $f$ yields $\vnorm{\nabla f(\v{x}_k)} \leq L$, from which we get \eqref{lemma1:MainIneq}.

\subsubsection*{Deduction of Theorem \ref{thm:CMD} from \eqref{lemma1:MainIneq}}
Let $\Xset^*$ be the set of optimizers of $f$ over $\Xset$, and $\v{x}^*$ be an arbitrary element in $\Xset^*$. Then, the  convexity of $f$ implies
\begin{align*}
\langle \nabla f(\v{x}_k), \v{x}^* - \v{x}_k \rangle \leq f^* - f(\v{x}_k).
\end{align*}
Using the above relation in \eqref{lemma1:MainIneq} with $\v{z} = \v{x}^*$ gives
\begin{align*}
D_\psi(\v{x}^*,\v{x}_{k+1} ) - D_\psi(\v{x}^*,\v{x}_k)  + \alpha_k [f(\v{x}_k)-f^*] \leq  \frac{\alpha_k^2 L^2}{2\mu}.
\end{align*} 
Summing the above over $k$ from 0 to $K$, we get
\begin{align*}
& D_\psi(\v{x}^*,\v{x}_{K+1} ) - D_\psi(\v{x}^*,\v{x}_0) \\
& \qquad +  \sum_{k=0}^K \alpha_k \left[ f(\v{x}_k)-f^* \right]
\leq \frac{L^2}{2\mu} \sum_{k=1}^K \alpha_k^2.
\end{align*} 
Taking $K \uparrow \infty$, the right hand side remains bounded, owing to the square summability of the $\alpha_k$'s. Bregman divergence is always nonnegative, and so is each summand in the second term  on the left hand side of the above inequality. Together, they imply that $\sum_{k=0}^\infty \alpha_k[ f(\v{x}_k)-f^* ] < \infty$ and 
$D_\psi(\v{x}_k, \v{x}^*)$ converges for each $\v{x}^* \in \Xset^*$. 
The non-summability of the $\alpha_k$'s further yields
\begin{align*}
\liminf_{k\rightarrow\infty} f(\v{x}_k) = f^*.
\end{align*}
Convergence of $D_\psi(\v{x}^*,\v{x}_k )$ for each $\v{x}^* \in \Xset^*$ implies the boundedness of the iterates $\v{x}_k$. Let $\{\v{x}_{k_\ell}\}_{\ell = 0}^\infty$ be the bounded subsequence of $\v{x}_{k}$'s along which
\begin{align}
\lim_{\ell\rightarrow\infty} f(\v{x}_{k_{\ell}}) =\liminf_{k\rightarrow\infty} f(\v{x}_k) = f^*.\label{thm:Eq1c}
\end{align}
This bounded sequence $\{\v{x}_{k_\ell}\}_{\ell = 0}^\infty$ has a (strongly) convergent subsequence. Function evaluations over that subsequence tend to $f^*$. Continuity of $f$ implies that the subsequence converges to a point in $\Xset^*$. Call this point $\v{x}^*$. Then, $D_\psi (\v{x}^*,\v{x}^{k} )$ converges, and it converges to zero over said subsequence, implying 
$$ \lim_{k \to \infty} D_\psi (\v{x}^*,\v{x}_{k} ) = 0.$$
Appealing to \eqref{eq:prop.bregman}, we conclude $\displaystyle \lim_{k \to \infty} \v{x}_k = \v{x}^*$.
This completes the proof of Theorem \ref{thm:CMD}.

\begin{remark}
\label{rem:InfDim}
This proof should be generalizable to $\Hset$ being infinite dimensional, where one would consider weak convergence of  $\v{x}_k$ to an optimizer of $f$ over $\Xset$. 
\end{remark}


\section{Distributed Mirror Descent}
\label{sec:DMD}
In this section, we consider a distributed variant of {\sf MD}. More precisely, we consider a collection of $N$ agents who collectively seek to minimize $f(\v{x}) := \sum_{i=1}^N f^i(\v{x})$ over $\Xset$. Agent $i$ only knows the convex but possibly non-smooth function $f^i$, and thus, the agents must solve the problem cooperatively.  The agents are allowed to exchange their iterates
only with their neighbors in an undirected graph $\Gfrak$.
Starting from $\v{x}_0^1, \ldots, \v{x}_0^N$, each agent communicates with its neighbors in $\Gfrak$ and updates its iterates $\v{x}^i_k$ at time $k$ as follows.
\begin{align}
\begin{aligned}
\v{v}^i_{k} &= \sum_{j=1}^N A_{ij} \v{x}^j_k, \\
\v{x}^i_{k+1} &= \argmin_{\v{z}\in\Xset}\left\{ \langle \nabla f(\v{v}^i_{k}), \v{z} - \v{v}^i_{k} \rangle+ \frac{1}{\alpha_k}D_\psi(\v{z},\v{v}^i_{k}) \right\}.
\end{aligned}
\label{eq:DMD}
\end{align}
Matrix $\v{A}$ thus encodes the communication graph $\Gfrak$, i.e.,
$A_{ij} \neq 0$ if and only if agent $j$ can communicate to agent
$i$ its current iterate, denoted by an edge between $i$ and $j$ in
$\Gfrak$. Rates for convergence of the function value in the above
\emph{distributed mirror descent} ({\sf DMD}) algorithm have been
reported in \cite[Theorem 2]{JLi2016}. We prove that {\sf DMD} drives
$\v{x}^1_k, \ldots \v{x}^N_k$ to a common optimizer $\v{x}^*$ of $f$
over $\Xset$.
\begin{theorem}
\label{thm:DMD} Suppose
\begin{itemize}
\item $f^i$ is $L$-Lipschitz continuous over $\Xset$,
\item $\{\alpha_k\}_{k=0}^\infty$ defines a nonincreasing sequence of positive step sizes that is square-summable, but not summable, i.e.,
$\sum_{k=0}^{\infty} \alpha_k = \infty,$ $\sum_{k=0}^{\infty} \alpha_k^2 < \infty$,
\item $\v{y} \mapsto D_\psi(\v{x}, \v{y})$ is convex,
\item $\v{A}$ is doubly stochastic, irreducible, and aperiodic.
\end{itemize}
Then, $\lim_{k \to \infty} \v{x}^i_k$ is identical across $i=1, \ldots, N$, and the limit optimizes $\sum_{i=1}^N f^i$ over $\Xset$ for $\v{x}_k$'s generated by {\sf DMD} in \eqref{eq:DMD}.
\end{theorem}
The first two assumptions are identical to the centralized counterpart in Section \ref{sec:centralized}. The third one is special to the distributed setting, and is crucial to the proof of the result. We remark that Bregman divergence $D_\psi$ is always strictly convex in its first argument. Our result requires convexity in the second argument. A sufficient condition is derived in \cite{Bauschke2001}, that requires $\psi$ to be thrice continuously differentiable and satisfy $\v{H}_\psi(\v{x}) \succeq 0$ and $\v{H}_\psi(\v{x}) + \nabla \v{H}_\psi(\v{x})(\v{x}-\v{y})\succeq 0$ for all $\v{x}$ and $\v{y}$ in $\Xset$, where $\v{H}_\psi$ stands for the Hessian of $\psi$. The last assumption defines a requirement on the information flow. Stated in terms of graph $\Gfrak$ that defines the connectivity among agents, it is sufficient to have $\Gfrak$ being connected with at least one node with a self-loop.


\subsection{Proof of Theorem \ref{thm:DMD}}
We first appeal to the optimality of $\v{x}^i_{k+1}$ in \eqref{eq:DMD} to conclude
\begin{align}
\langle \alpha_k \nabla f^i(\v{v}^i_k) + \nabla_1 D_\psi(\v{x}^i_{k+1},\v{v}^i_k), \v{z}-\v{x}^i_{k+1} \rangle \geq 0
\label{eq:optDMD}
\end{align}
for every $\v{z} \in \Xset$. The properties of Bregman divergence in \eqref{eq:prop.bregman} yield
\begin{align}
\begin{aligned}
& \langle \nabla D_\psi( \v{x}^i_{k+1},\v{v}^i_k), \v{z} - \v{x}^i_{k+1} \rangle \\
& \qquad = \langle \nabla \psi(\v{x}^i_{k+1}) - \nabla \psi(\v{v}_k), \v{z} - \v{x}^i_{k+1} \rangle \\
& \qquad = D_\psi(\v{z}, \v{v}^i_{k}) - D_\psi( \v{z},\v{x}^i_{k+1}) - D_\psi(\v{v}^i_{k}, \v{x}^i_{k+1}).
\end{aligned}
\label{eq:nablaD}
\end{align}
Substituting the above equality in \eqref{eq:optDMD}, and summing over $i=1, \ldots, N$, we get an inequality of the form
\begin{align}
\sum_{i=1}^N S^i_k (\v{z}) + \sum_{i=1}^N T^i_k(\v{z}) \geq 0,
\label{eq:SikzTikz}
\end{align}
where
\begin{align*}
S^i_k (\v{z}) &:= \langle \alpha_k \nabla f^i(\v{v}^i_k), \v{z} - \v{x}^i_{k+1} \rangle, \\
T^i_k(\v{z}) &:= D_\psi(\v{z}, \v{v}^i_{k}) - D_\psi( \v{z},\v{x}^i_{k+1}) - D_\psi(\v{x}^i_{k+1}, \v{v}^i_{k})
\end{align*}
for each $\v{z} \in \Xset$. We provide upper bounds on each of the summations separately. In the sequel, we use the notation
\begin{align}
\ol{\v{x}}_k = \frac{1}{N} \sum_{i=1}^N \v{x}^i_k.
\end{align}


\subsubsection*{An upper bound for $\sum_{i=1}^N S^i_k (\v{z})$}
\begin{align*}
\sum_{i=1}^N S^i_k (\v{z})
&= \sum_{i=1}^N \langle \alpha_k \nabla f^i(\v{v}^i_k), \v{z} - \v{v}^i_k \rangle \\
& \qquad + \sum_{i=1}^N \langle \alpha_k \nabla f^i(\v{v}^i_k), \v{v}^i_k - \v{x}^i_{k+1} \rangle.
\end{align*}
We bound each summand in both the summations above. Using the convexity of $f^i$, we get
\begin{align}
& \langle \nabla f^i(\v{v}^i_k), \v{z} - \v{v}^i_k \rangle  \notag\\
& \qquad \leq f^i(\v{z}) - f^i(\v{v}^i_{k}) \notag \\
& \qquad = f^i(\v{z}) - f^i(\ol{\v{x}}_{k}) + f^i(\ol{\v{x}}_{k}) - f^i(\v{v}^i_k) \notag \\
& \qquad \leq f^i(\v{z}) - f^i(\ol{\v{x}}_{k}) + \langle \nabla f^i(\v{v}^i_k), \ol{\v{x}}_{k} - \v{v}^i_k \rangle \notag \\
& \qquad \leq f^i(\v{z}) - f^i(\ol{\v{x}}_{k}) + L \vnorm{\ol{\v{x}}_{k} - \v{v}^i_k}.
\label{eq:Sikz.1}
\end{align}
The last line follows from the Cauchy-Schwarz inequality and that $f^i$ is $L$-Lipschitz. Further,
the Cauchy-Schwarz and arithmetic-geometric mean inequalities yield
\begin{align}
&\langle \alpha_k \nabla f^i(\v{v}^i_k), \v{v}^i_k - \v{x}^i_{k+1} \rangle \notag \\
& \qquad \leq \frac{\alpha_{k}^2}{2\mu}{\vnorm{\nabla f^i(\v{v}^i_k)}^2} + \frac{\mu}{2} \vnorm{\v{v}^i_k - \v{x}^i_{k+1}}^2 \notag\\
& \qquad \leq \frac{\alpha_{k}^2 L^2}{2\mu} + \frac{\mu}{2} \vnorm{\v{v}^i_k - \v{x}^i_{k+1}}^2.
\label{eq:Sikz.2}
\end{align}
The last line is a consequence of $f^i$ being $L$-Lipschitz.

Combining \eqref{eq:Sikz.1} and \eqref{eq:Sikz.2} then allows us to deduce
\begin{align}
\sum_{i=1}^N S^i_k (\v{z})
& \leq \alpha_k [f(\v{z}) - f(\ol{\v{x}}_{k})] + \alpha_k L \sum_{i=1}^N \vnorm{\ol{\v{x}}_{k} - \v{v}^i_k} \notag \\
& \qquad + \alpha_{k}^2  \frac{N L^2}{2\mu} + \frac{\mu}{2} \sum_{i=1}^N \vnorm{\v{v}^i_k - \v{x}^i_{k+1}}^2.
\label{eq:Sikz.3}
\end{align}
$\v{A}$ is doubly stochastic and $\vnorm{\cdot}$ is a convex function, implying
$$ \sum_{i=1}^N \vnorm{\ol{\v{x}}_{k} - \v{v}^i_k} = \sum_{i=1}^N \vnorm{\ol{\v{x}}_{k} - \sum_{j=1}^N A_{ij}\v{x}^j_k} \leq \sum_{i=1}^N \vnorm{\ol{\v{x}}_{k} - \v{x}^i_k}.$$
Upon utilizing the above in \eqref{eq:Sikz.3}, we derive the required upper bound on $\sum_{i=1}^N S^i_k (\v{z})$.
\begin{align}
\sum_{i=1}^N S^i_k (\v{z})
& \leq \alpha_k [f(\v{z}) - f(\ol{\v{x}}_{k})] + \alpha_k L \sum_{i=1}^N \vnorm{\ol{\v{x}}_{k} - \v{x}^i_k} \notag \\
& \qquad + \alpha_{k}^2  \frac{N L^2}{2\mu} + \frac{\mu}{2} \sum_{i=1}^N \vnorm{\v{v}^i_k - \v{x}^i_{k+1}}^2.
\label{eq:Sikz.4}
\end{align}


\subsubsection*{An upper bound for $\sum_{i=1}^N T^i_k (\v{z})$}
\begin{align*}
\sum_{i=1}^N T^i_k(\v{z})
& = \sum_{i=1}^N \left[ D_\psi(\v{z}, \v{v}^i_{k}) - D_\psi( \v{z},\v{x}^i_{k+1}) \right] \\
& \qquad - \sum_{i=1}^N D_\psi(\v{x}^i_{k+1}, \v{v}^i_{k}).
\end{align*}
To bound the right-hand side, we utilize the convexity of $D_\psi$ in the second argument and the doubly stochastic nature of $\v{A}$ to obtain
\begin{align*}
\sum_{i=1}^N D_\psi(\v{z}, \v{v}^i_{k})
= \sum_{i=1}^N D_\psi\left(\v{z}, \sum_{j=1}^N A_{ij}\v{x}^i_{k}\right) \leq \sum_{i=1}^N D_\psi\left(\v{z}, \v{x}^i_{k}\right).
\end{align*}
Also, \eqref{eq:prop.bregman} gives
\begin{align*}
D_\psi(\v{x}^i_{k+1}, \v{v}^i_{k}) \geq \frac{\mu}{2} \vnorm{\v{v}^i_k - \v{x}^i_{k+1}}^2.
\end{align*}
Combination of the above two inequalities then provides the sought upper bound on $\sum_{i=1}^N T^i_k(\v{z})$.
\begin{align}
\sum_{i=1}^N T^i_k(\v{z})
& \leq \sum_{i=1}^N \left[ D_\psi(\v{z}, \v{x}^i_{k}) - D_\psi( \v{z},\v{x}^i_{k+1}) \right] \notag \\
& \qquad - \frac{\mu}{2} \sum_{i=1}^N   \vnorm{\v{v}^i_k - \v{x}^i_{k+1}}^2.
\label{eq:Tikz.3}
\end{align}


\subsubsection*{Utilizing the bounds on $\sum_{i=1}^N S^i_k (\v{z})$ and $\sum_{i=1}^N T^i_k (\v{z})$}
Let $\v{x}^* \in \Xset^*$ be an arbitrary optimizer of $f$ over $\Xset$.
Applying the bounds in \eqref{eq:Sikz.4} and \eqref{eq:Tikz.3} in \eqref{eq:SikzTikz} with $\v{z} = \v{x}^*$ then gives
\begin{align*}
&\sum_{i=1}^N \left[ D_\psi( \v{x}^*,\v{x}^i_{k+1}) - D_\psi(\v{x}^*, \v{x}^i_{k}) \right] + \alpha_k [f(\ol{\v{x}}_{k}) - f(\v{x}^*)]  \\
& \qquad \leq  \alpha_k L \sum_{i=1}^N \vnorm{\ol{\v{x}}_{k} - \v{x}^i_k} + \alpha_k^2 \frac{N L^2}{2\mu}.
\end{align*}
Summing the above over $k = 0, \ldots, K$, we obtain
\begin{align}
&\sum_{i=1}^N \left[ D_\psi( \v{x}^*,\v{x}^i_{K+1}) - D_\psi(\v{x}^*, \v{x}^i_{0}) \right] + \sum_{k=1}^K \alpha_k [f(\ol{\v{x}}_{k}) - f(\v{x}^*)]  \notag \\
& \qquad \leq  L \sum_{k=0}^K \sum_{i=1}^N \alpha_k  \vnorm{\ol{\v{x}}_{k} - \v{x}^i_k} + \frac{N L^2}{2\mu} \sum_{k=0}^K \alpha_k^2.
\label{eq:combined}
\end{align}
We mimic the style of arguments in the proof of Theorem \ref{thm:CMD} to complete the derivation, and provide an upper bound on the double summation on the right hand side of the above inequality.


\subsubsection*{An upper bound on $\sum_{k=0}^K \sum_{i=1}^N \alpha_k \vnorm{\ol{\v{x}}_{k} - \v{x}^i_k}$}

To derive this bound, we fix an orthonormal basis for $\Hset$ with $\dim \Hset = d$. Let $\mathpzc{x}^i_k \in \Rset^d$ denote the coordinates of $\v{x}^i_k$ in that basis. The coordinates for the centroid $\ol{\v{x}}_k$ are given by $\ol{\mathpzc{x}}_k$. The inner product in $\Hset$ becomes the usual dot product among the corresponding coordinates. The norm becomes the usual Euclidean 2-norm in the coordinates. Define
$$\v{X}_k^\T := \begin{pmatrix} \mathpzc{x}^1_k  \ldots \mathpzc{x}^N_k \end{pmatrix} \in \Rset^{d \times N}.$$
Let $\bone \in \Rset^{N}$ denote a vector of all ones and $\v{I}\in\Rset^{N \times N}$ be the identity matrix. Also, define
$$ \v{P} := \v{I} - \frac{1}{N}\bone \bone^\T $$
for convenience. Equipped with this notation, we then have
\begin{align}
\sum_{i=1}^N \vnorm{\ol{\v{x}}_{k} - \v{x}^i_k}
& = \sum_{i=1}^N \vnorm{\ol{\mathpzc{x}}_{k} - \mathpzc{x}^i_k}_2 \notag \\
& \leq \sqrt{N} \vnorm{\v{P} \v{X}_k}_F  \notag\\
& \leq  \underbrace{\sqrt{N \cdot\min\{N, d\}}}_{:=N'} \vnorm{\v{P} \v{X}_k}_2. \label{eq:norm2}
\end{align}
Here, $\vnorm{\cdot}_F$ and $\vnorm{\cdot}$ denote the Frobenius and the 2-norm of matrices, respectively. In what follows, we bound $ \sum_{k=0}^K \alpha_k \vnorm{\v{P}\v{X}_k}_2 $ from above.

Collecting the coordinates of $\v{v}^i_k$ in $\v{V}_k \in \Rset^{N \times d}$ similarly to $\v{X}_k$, we have
\begin{align}
\vnorm{\v{P} \v{X}_{k+1}}_2
&= \vnorm{ \v{P} \left( \v{A} \v{X}_{k}  + \v{X}_{k+1} - \v{V}_k \right)}_2 \notag \\
& \; \leq \vnorm{ \v{A} \v{P} \v{X}_{k} }_2 + \vnorm{ \v{P} \left(\v{X}_{k+1} - \v{V}_k\right)}_2
\label{eq:Xk1.1}
\end{align}
since $\v{A}$ commutes with $\v{P}$. We bound each term on the right hand side above.
For the first term, notice that
$$\v{P} \v{X}_k = \left( \v{I} - \frac{1}{N}\bone \bone^\T \right) \v{X}_{k} \perp \bone.$$
Since $\v{A}$ is doubly stochastic (i.e., $\v{A} \bone = \bone$), the Perron-Frobenius theorem \cite[Theorem 8.4.4]{HJbook} and the Courant-Fischer theorem \cite[Theorem 4.2.11]{HJbook} together yield
\begin{align}
\vnorm{ \v{A} \v{P} \v{X}_{k} }_2
 \leq  \sigma_2(\v{A}) \vnorm{\v{P} \v{X}_{k} }_2  \label{eq:Xk1.2},
\end{align}
where $\sigma_2(\v{A})$ is the second largest singular value of $\v{A}$. Further, $\v{A}$ is irreducible, and aperiodic, implying $\sigma_2(\v{A}) < 1$.
To bound the second term on the right hand side of \eqref{eq:Xk1.1}, we use that matrix norms are submultiplicative, and hence we have
\begin{align}
\vnorm{ \v{P} \left(\v{X}_{k+1} - \v{V}_k\right)}_2
& \leq {\vnorm{\v{P}}_2} \vnorm{\v{X}_{k+1} - \v{V}_k}_2  \notag \\
& = \sum_{i=1}^N  \vnorm{\v{v}^i_k - \v{x}^i_{k+1}},
\label{eq:VX}
\end{align}
because $\vnorm{\v{P}}_2 = 1$.
To bound each term in the above summation, we utilize \eqref{eq:optDMD} and \eqref{eq:nablaD} with $\v{z} = \v{v}^i_k$ to obtain
\begin{align}
& \langle \alpha_k \nabla f^i(\v{v}^i_k), \v{v}^i_k-\v{x}^i_{k+1} \rangle \notag\\
& \qquad \geq \langle \nabla \psi(\v{v}_k) - \nabla \psi(\v{x}^i_{k+1}), \v{v}^i_k - \v{x}^i_{k+1} \rangle.
\label{eq:VX.1}
\end{align}
Since $f^i$ is $L$-Lipschitz, and $\psi$ is $\mu$-strongly convex, we have the following two inequalities
\begin{align*}
{\langle \alpha_k \nabla f^i(\v{v}^i_k), \v{v}^i_k-\v{x}^i_{k+1} \rangle}
&\leq {\alpha_k L \vnorm{\v{v}^i_k-\v{x}^i_{k+1}}}, \\
{\langle \nabla \psi(\v{v}_k) - \nabla \psi(\v{x}^i_{k+1}), \v{v}^i_k - \v{x}^i_{k+1} \rangle}
&\geq { \mu \vnorm{\v{v}^i_k - \v{x}^i_{k+1}}^2}.
\end{align*}
that together with \eqref{eq:VX.1} gives
\begin{align*}
 \vnorm{\v{v}^i_k - \v{x}^i_{k+1}} \leq \frac{\alpha_k L}{\mu}.
\end{align*}
Summing the above over $i=1,\ldots,N$, we obtain an upper bound on the right hand side of \eqref{eq:VX}. Utilizing that bound and \eqref{eq:Xk1.2} in \eqref{eq:Xk1.1}, we get
\begin{align}
\vnorm{\v{P}\v{X}_{k+1}}_2
\leq \sigma_2(\v{A}) \vnorm{ \v{P} \v{X}_{k}}_2
+ \alpha_k \frac{NL}{\mu}.
\end{align}

For convenience we suppress the dependency of $\sigma_2$ on $\v{A}$ in the sequel. Iterating the above inequality gives
\begin{align*}
\vnorm{\v{P}\v{X}_{k}}_2
\leq \sigma_2^k \vnorm{ \v{P} \v{X}_{0}}_2
+  \frac{NL}{\mu} \sum_{\ell=0}^{k-1} \alpha_\ell \sigma_2^{k-\ell-1}
\end{align*}
for each $k \geq 1$,
which further yields
\begin{align}
\sum_{k=0}^K \alpha_k \vnorm{\v{P}\v{X}_{k} }_2
& \leq \alpha_0 \vnorm{  \v{P} \v{X}_{0} }_2 + \sum_{k=1}^K \alpha_k \sigma_2^k \vnorm{  \v{P} \v{X}_{0}}_2 \notag \\
& \qquad +  \frac{NL}{\mu} \sum_{k=1}^K \alpha_k \sum_{\ell=0}^{k-1} \alpha_\ell \sigma_2^{k-\ell-1}.
\label{eq:XkEk.1}
\end{align}
Now, $\sigma_2 < 1$ and the $\alpha_k$'s are nonincreasing. Thus, we have
$$\alpha_0 + \sum_{k=1}^K  \alpha_k \sigma_2^k \leq \alpha_0 \sum_{k=0}^\infty \sigma_2^k = \alpha_0 (1-\sigma_2)^{-1},$$
and using this in \eqref{eq:XkEk.1} gives the following required bound.
\begin{align}
& \sum_{k=0}^K \sum_{i=1}^N \alpha_k \vnorm{\ol{\v{x}}_{k} - \v{x}^i_k} \notag \\
& \quad \leq N'\sum_{k=0}^K \alpha_k \vnorm{\v{P} \v{X}_{k}}_2 \notag \\
&\quad \leq \frac{\alpha_0 N'}{1-\sigma_2} \vnorm{ \v{P} \v{X}_{0}}_2 +  \frac{NN'L}{\mu} \sum_{k=1}^K \sum_{\ell=0}^{k-1} \alpha_\ell^2 \sigma_2^{k-\ell-1} \notag \\
&\quad \leq \frac{\alpha_0 N'}{1-\sigma_2} \vnorm{ \v{P} \v{X}_{0}}_2 +  \frac{NN'L}{\mu} \sum_{\ell=0}^{K-1} \alpha_\ell^2 \sum_{k=\ell + 1}^K \sigma_2^{k-\ell-1} \notag \\
&\quad \leq \frac{\alpha_0 N'}{1-\sigma_2} \vnorm{  \v{P} \v{X}_{0} }_2 +  \frac{NN'L}{\mu (1 - \sigma_2)} \sum_{k=0}^{K-1} \alpha_k^2. \label{eq:XkEk.2}
\end{align}

\subsubsection*{Proof of the result by combining all upper bounds}
Utilizing \eqref{eq:XkEk.2} in \eqref{eq:combined} yields
\begin{align}
&\sum_{i=1}^N \left[ D_\psi( \v{x}^*,\v{x}^i_{K+1}) - D_\psi(\v{x}^*, \v{x}^i_{0}) \right] + \sum_{k=0}^K \alpha_k [f(\ol{\v{x}}_{k}) - f(\v{x}^*)]  \notag \\
&\qquad \leq \frac{ \alpha_0 N' L}{1-\sigma_2} \vnorm{  \v{P} \v{X}_{0}}_2 +  \frac{N N' L^2}{\mu (1 - \sigma_2)} \sum_{k=0}^{K-1} \alpha_k^2
+ \frac{N L^2}{2\mu} \sum_{k=0}^K \alpha_k^2. \label{eq:combined.2}
\end{align}
Driving $K\uparrow \infty$, the right hand sides of \eqref{eq:XkEk.2} and \eqref{eq:combined.2} converge as the sequence of $\alpha_k$'s are square summable. Further, the $\alpha_k$'s are non-summable, and hence, we conclude
\begin{enumerate}[leftmargin=*]
\item $\liminf_{K \to \infty} \sum_{i=1}^N  \vnorm{\ol{\v{x}}_{K} - \v{x}^i_K} = 0$,
\item $ \liminf_{K \to \infty} f(\ol{\v{x}}_K) = f^*$,
\item $\sum_{i=1}^N  D_\psi( \v{x}^*,\v{x}^i_{K})$ converges.
\end{enumerate}
Convergence of the Bregman divergences for each $\v{x}^* \in \Xset^*$ implies the boundedness of the iterates $\v{x}^i_K$ for each $i=1, \ldots,N$.
Recall that $\v{X}_K$ denotes the collective iterate for all agents at time $K$.
Consider the bounded subsequence of $\v{X}_{K}$, denoted $\v{X}_{K_\ell}$, along which
\begin{align*}
\lim_{\ell\to\infty} \sum_{i=1}^N  \vnorm{\ol{\v{x}}_{K_\ell} - \v{x}^i_{K_\ell}}
&= \liminf_{K\to\infty} \sum_{i=1}^N  \vnorm{\ol{\v{x}}_{K} - \v{x}^i_K} = 0, \\
\lim_{\ell\to\infty} f(\ol{\v{x}}_{K_{\ell}})
&= \liminf_{K\rightarrow\infty} f(\ol{\v{x}}_K) = f^*.
\end{align*}
This bounded sequence $\{\v{X}_{K_\ell}\}_{\ell = 0}^\infty$ has a (strongly) convergent subsequence. Over that subsequence, the agents' iterates converge to the centroid, and the function evaluations over the centroid tend to $f^*$. Continuity of $f$ implies that each agent's iterate over that subsequence converges to the \textit{same} point in $\Xset^*$. Call this point $\v{x}^*$. Since $D_\psi (\v{x}^*,\v{x}^i_{K} )$ converges, it converges to zero over that subsequence, implying
$ \lim_{k \to \infty} D_\psi (\v{x}^*,\v{x}^i_{K} ) = 0$
for each $i$. Appealing to \eqref{eq:prop.bregman}, we infer $\displaystyle \lim_{k \to \infty} \v{x}^i_K = \v{x}^*$, concluding the proof.


\section{Numerical experiments}
\label{sec:eg}
Theorems \ref{thm:CMD} and \ref{thm:DMD} guarantee the convergence of the iterates to the optimizer, but do not provide convergence rates with non-summable but square summable step-sizes. Given the lack of rates, we empirically illustrate that mirror descent -- both in centralized and distributed settings -- often outperforms vanilla subgradient methods on simple examples with our step sizes. Our simulations are different from many prior works, e.g., \cite{Nedic2014}, where they choose $\alpha_k = \frac{a}{\sqrt{k+1}}, a>0$ to guarantee the fastest convergence of the function value.

\begin{figure}[H]
\centering
    \begin{subfigure}[b]{0.5\textwidth}
        \centering
    \includegraphics[width=.42\textwidth]{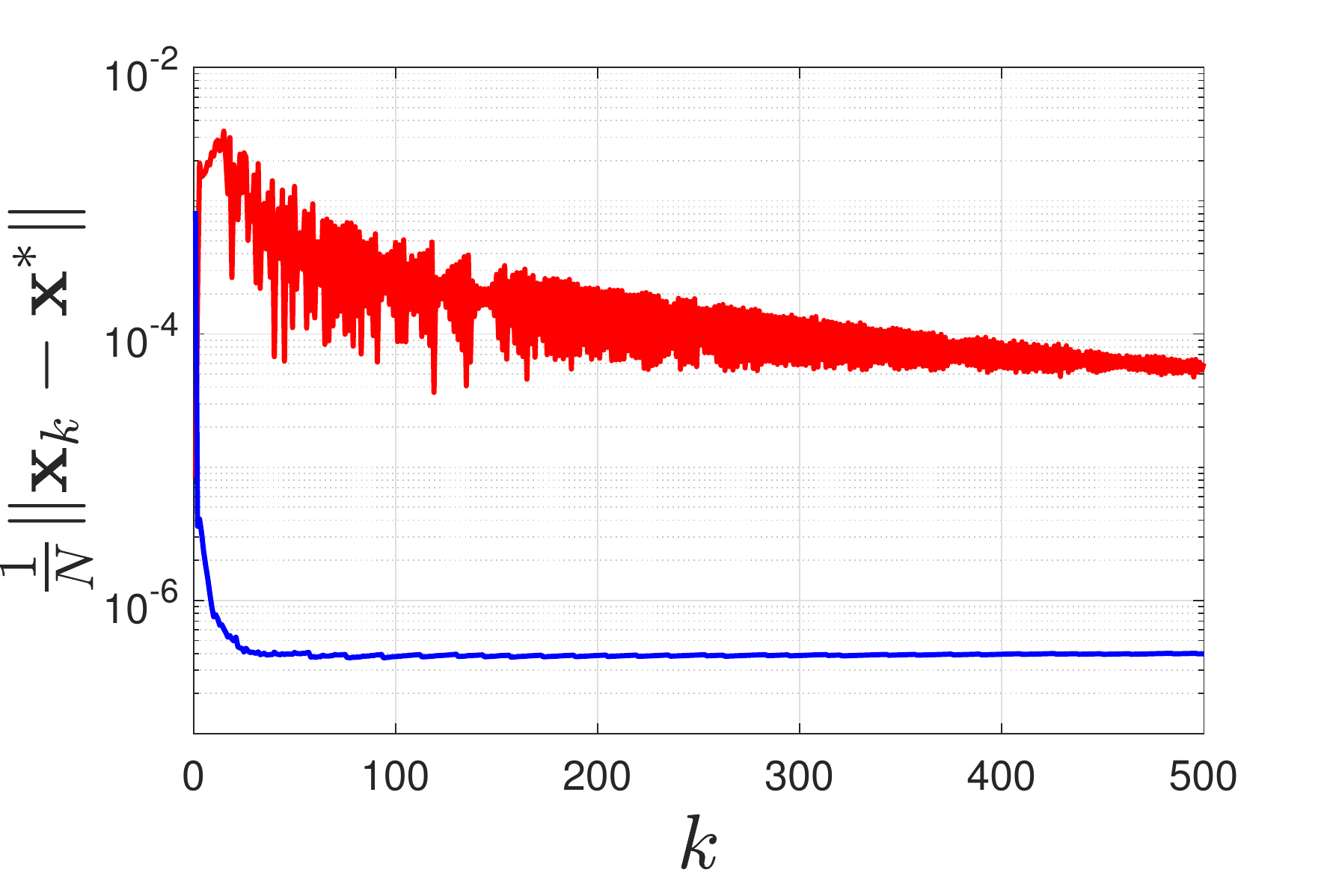}
    \caption{Centralized algorithms.}
    \label{fig:CMD}
    \end{subfigure}
    \begin{subfigure}[b]{0.5\textwidth}
    \centering
        \includegraphics[width=0.50\textwidth]{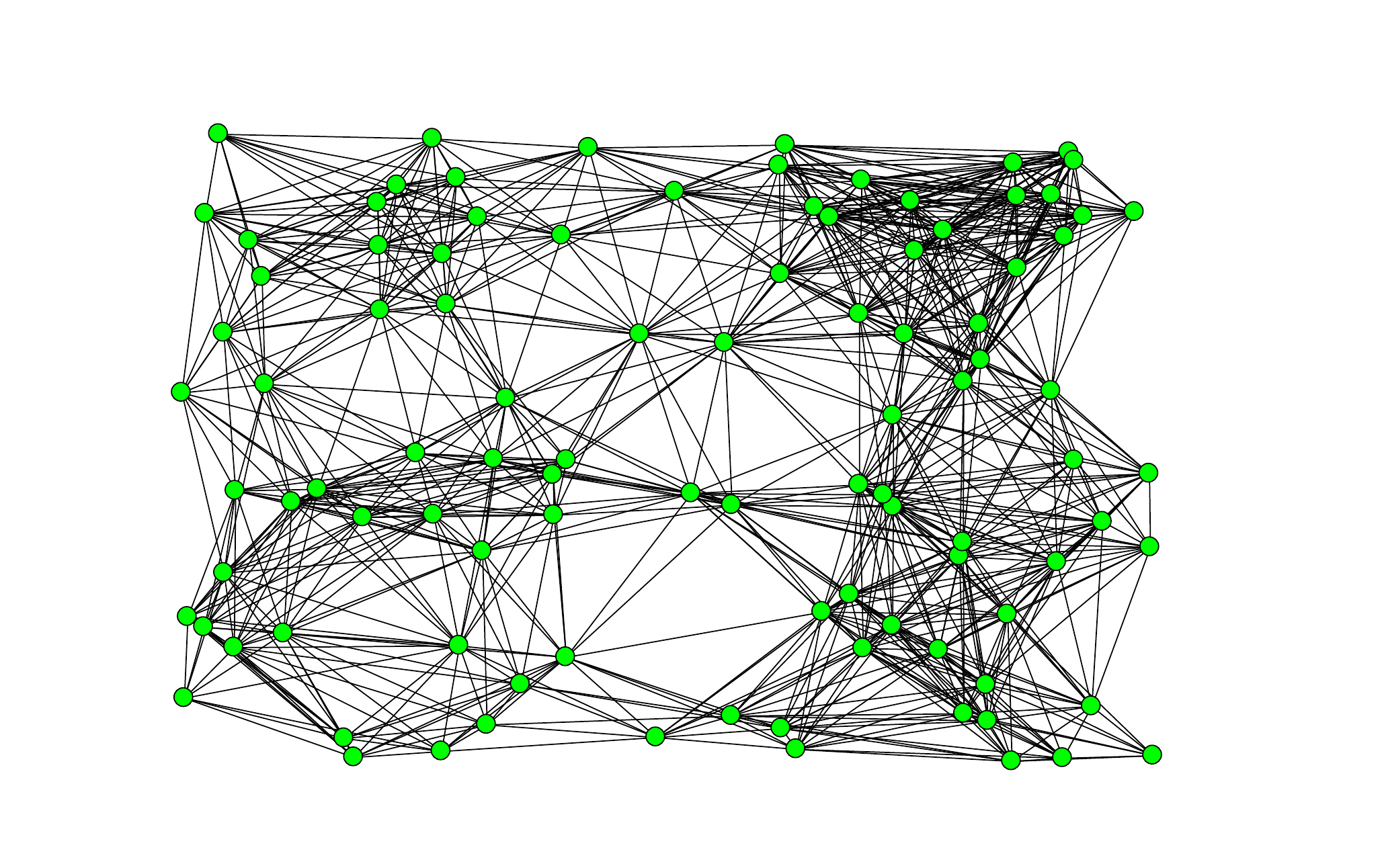}
        \includegraphics[width=0.42\textwidth]{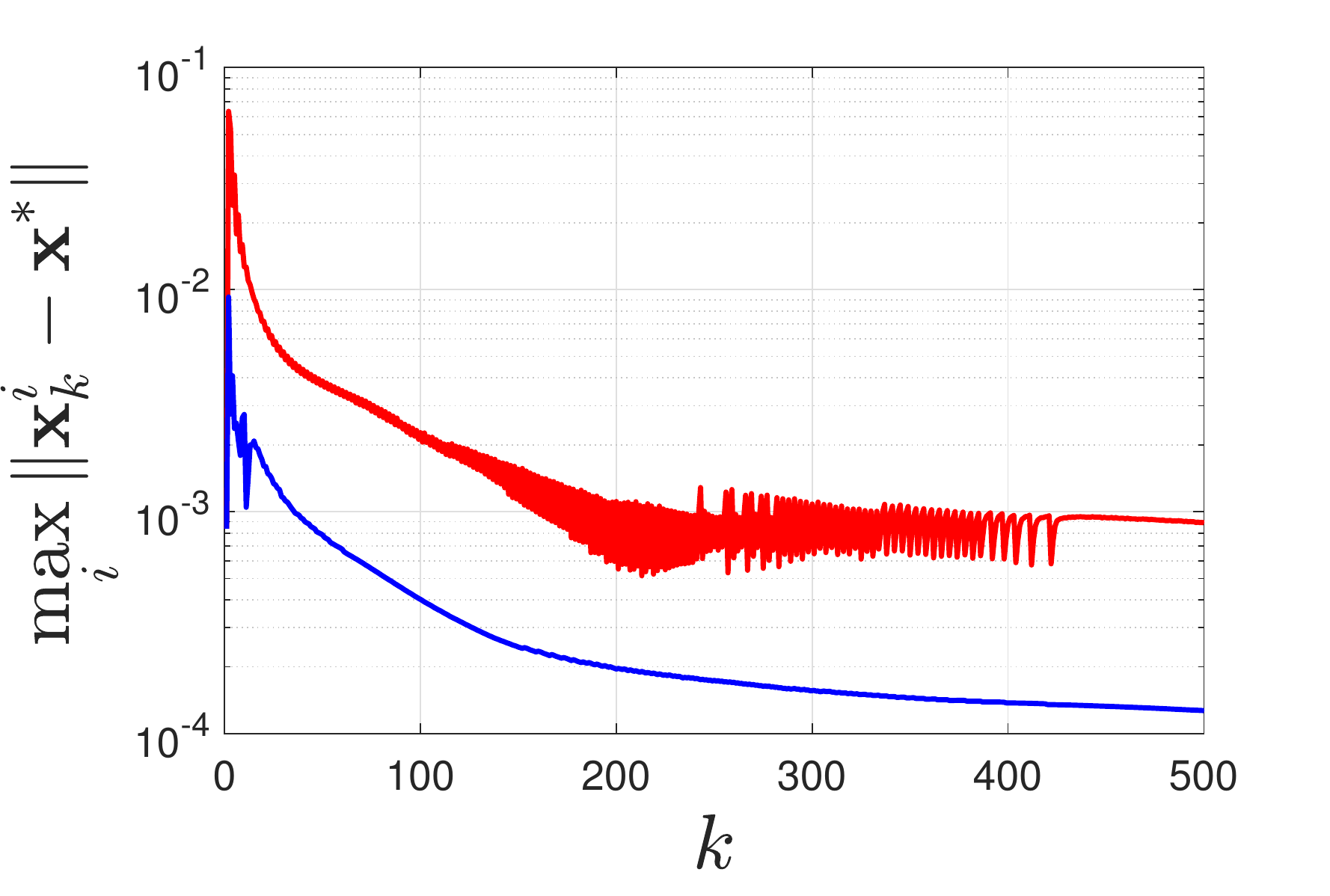}
        \caption{Distributed algorithms on a network with 939 edges.}
        \label{fig:network_r03}
    \end{subfigure}
    \begin{subfigure}[b]{0.5\textwidth}
        \centering
        \includegraphics[width=0.50\textwidth]{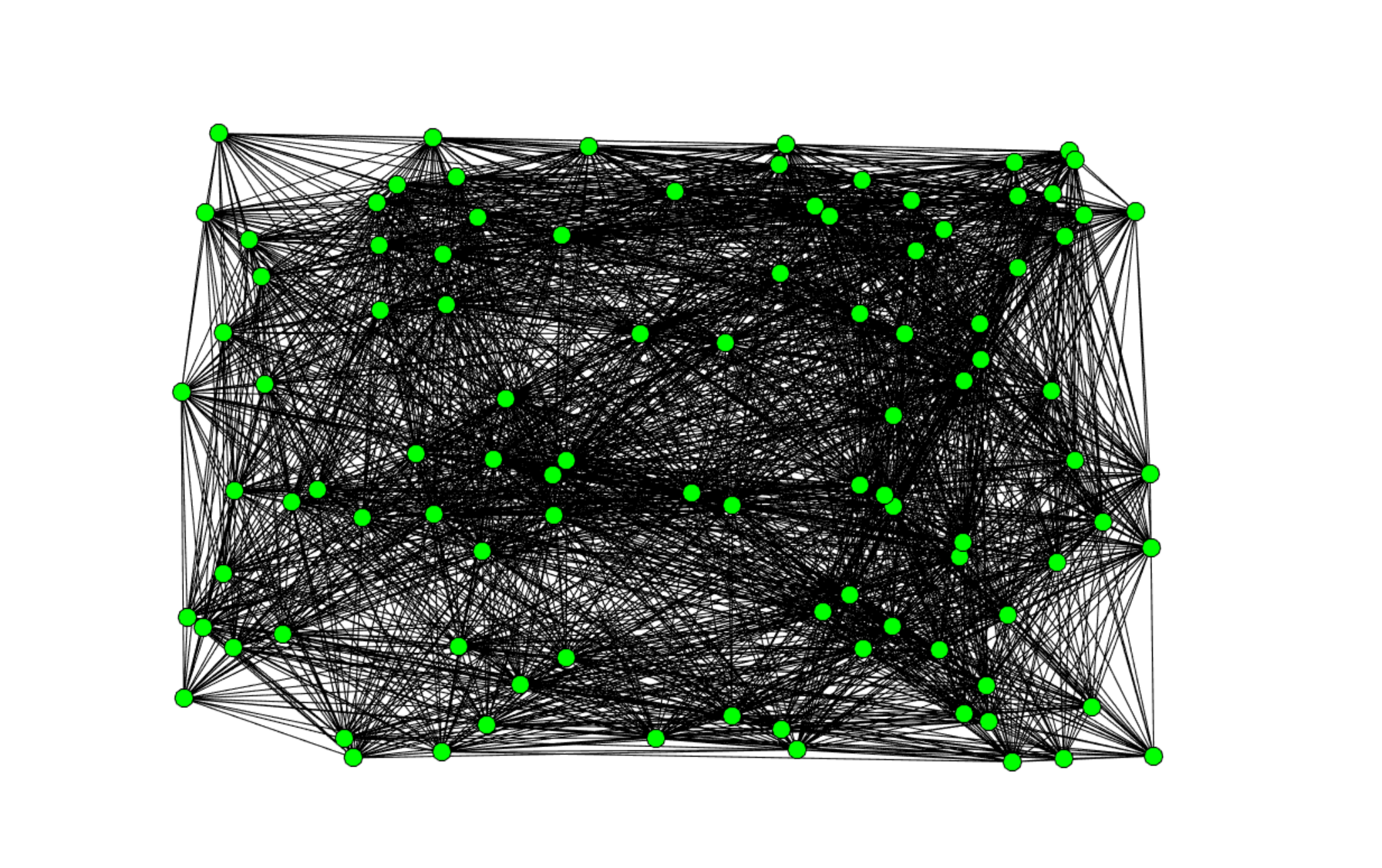}
        \includegraphics[width=0.42\textwidth]{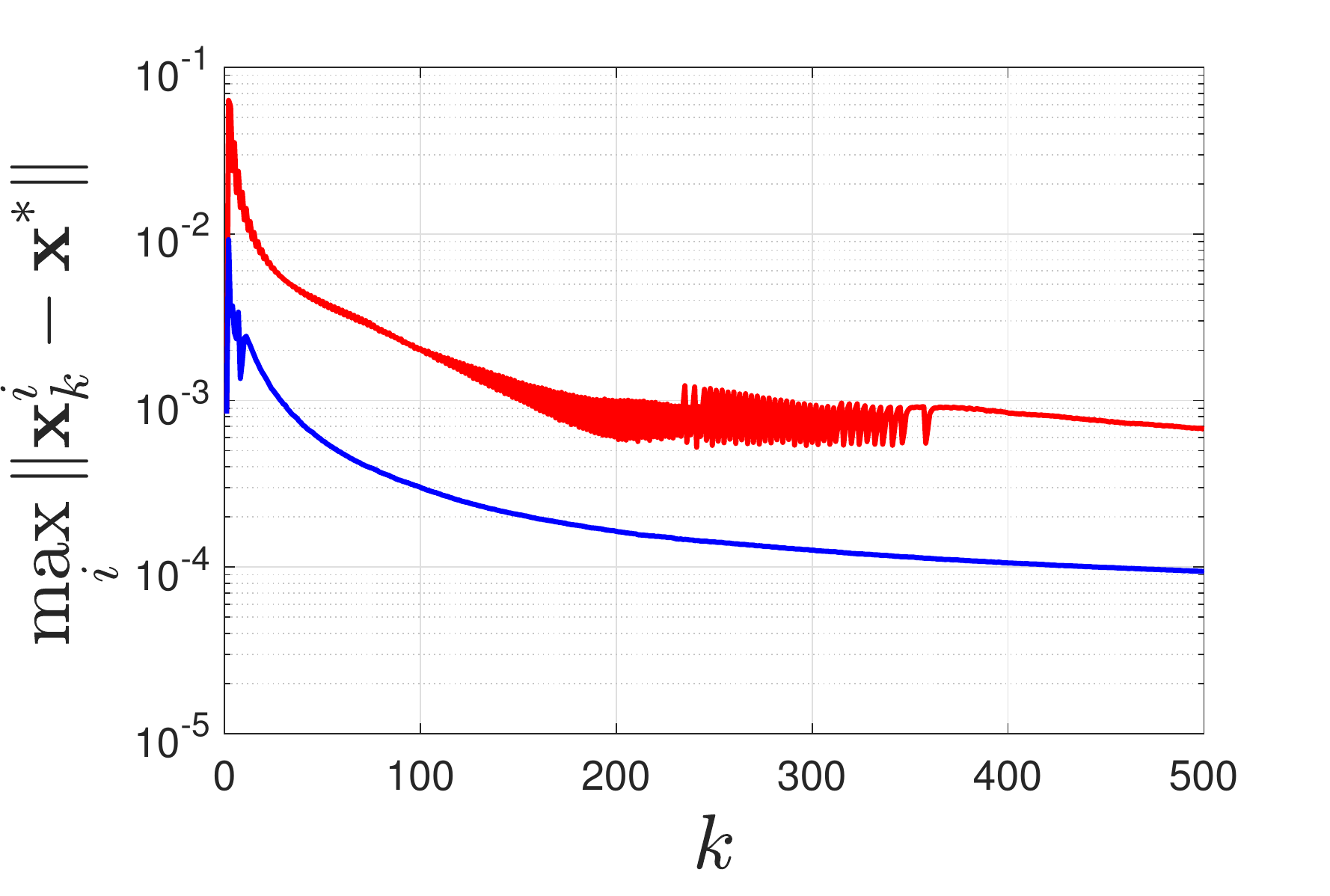}
        \caption{Distributed algorithms on a network with 2678 edges.}
        \label{fig:network_r06}
    \end{subfigure}
    \caption{Diagram showing the convergence behavior of projected subgradient method (\protect\redline) and entropic mirror descent (\protect\blueline) for \eqref{eq:robReg} with $N=100$ and $d=10$. All entries of $\v{G}$ and $\v{h}$ are chosen uniformly at random from $[0,1]$. The algorithm was initialized at a random point in $\Xset$. Plot (a) shows the dynamics of the centralized algorithms. Plots (b) and (c) show dynamics of the distributed variants over the respective networks with 939 and 2678 edges, respectively.}
\label{fig:MD}
\end{figure}

Consider the following robust linear regression problem over a simplex.
\begin{align}
\underset{\v{x} \in \Rset^d}{\text{minimize}}   \ \
 \vnorm{\v{G}\v{x} - \v{h}}_1, \ \ \text{subject to} \ \ \bone^\T \v{x} = 1, \ \v{x} \geq 0.
\label{eq:robReg}
\end{align}
Robust regression fits a linear model to the data $\v{G} \in \Rset^{N \times d}, \v{h} \in \Rset^N$. It differs from ordinary least squares in that the objective function penalizes the entry-wise absolute deviation from the linear fit rather than the squared residue, and is known to be robust to outliers \cite{karst1958linear}. Consider two different Bregman divergences on the $d$-dimensional simplex $\Xset$ defined by the
Euclidean distance $\psi_1(\v{x}) := \frac{1}{2}\vnorm{\v{x}}_2^2$, and
negative entropy $\psi_2(\v{x}) := \sum_{j=1}^d x^j \log x^j$.
Centralized mirror descent with $D_{\psi_1}$ amounts to a projected subgradient algorithm where each iteration is a subgradient step followed by a projection on $\Xset$. With $D_{\psi_2}$, the updates define an exponentiated gradient method, also known as the entropic mirror descent algorithm (cf. \cite{precup1997exponentiated, BeckT2003}). Its updates are given by
$$ {x}^j_{k+1} := \frac{x^j_k \exp\left( - \alpha_k [\nabla f(\v{x}^k)]^j \right)}{\sum_{\ell=1}^{d} x^\ell_k {\exp\left( - \alpha_k [\nabla f(\v{x}^k)]^\ell \right)}},$$
where the objective in \eqref{eq:robReg} is $f(\v{x})$, and
$$ \nabla f(\v{x}) = \sum_{j=1}^N \sgn([\v{g}^{i}]^\T \v{x} - h^{i}) \v{g}^{i}.$$
Here, $\sgn(\cdot)$ denotes the sign of the argument, and $[\v{g}^i]^\T$ is the $i$-th row of $\v{G}$.
Negative entropy being a `natural' function over simplex, entropic
mirror descent enjoys faster convergence than projected subgradient
descent, as shown in Figure \ref{fig:CMD} using step sizes $\alpha_k
= \frac{1}{5(k+1)}$.

Next, consider the case where each node $i=1, \ldots, N$ in a graph
only knows $\v{g}^i$ and $h^i$, and they together seek to minimize
$\sum_{i=1}^N \lvert [\v{g}^i]^\T \v{x} - h^i \rvert$. Figures
\ref{fig:network_r03} and \ref{fig:network_r06} show how the
distributed variant of entropic mirror descent outperforms that of
projected subgradient method with steps-sizes $\alpha_k =
\frac{1}{5(k+1)}$. We choose $\v{A}$ as the transition probabilities
of a Markov chain in the Metropolis-Hastings algorithm for the
respective graphs. Centralized algorithms converge faster than
distributed algorithms; however, the denser the graph, the faster
the convergence is of the distributed algorithms.
%


\section{Conclusion}
\label{sec:conclusion}
In this paper, we proved guaranteed convergence of the iterates (i.e., the problem variables) to an optimizer in {\sf MD} on a finite dimensional Hilbert space, using a specific choice of step size in both centralized and distributed settings. The convergence holds even when minimizing possibly non-smooth and non-strongly convex functions. This convergent behavior generalizes a similar property of subgradient methods. Extension to the case with additive noise with bounded support, and to infinite dimensional Hilbert and Banach spaces remain interesting directions for future research.

\bibliographystyle{plain}
\bibliography{refs}

\end{document}